\documentclass[11pt]{amsart}

\usepackage{pb-diagram}

\newtheorem{thm}{Theorem}[section]   
\newtheorem{cor}[thm]{Corollary}     
\newtheorem{lem}[thm]{Lemma}         
\newtheorem{prop}[thm]{Proposition}  

\theoremstyle{definition}
\newtheorem{defn}[thm]{Definition}   

\theoremstyle{remark}
\newtheorem{rem}[thm]{Remark}        

\numberwithin{equation}{section}     

\newcommand{\secref}[1]{Section~\textup{\ref{#1}}}
\newcommand{\thmref}[1]{Theorem~\textup{\ref{#1}}}
\newcommand{\corref}[1]{Corollary~\textup{\ref{#1}}}
\newcommand{\lemref}[1]{Lemma~\textup{\ref{#1}}}
\newcommand{\propref}[1]{Proposition~\textup{\ref{#1}}}
\newcommand{\defnref}[1]{Definition~\textup{\ref{#1}}}

\newcommand{\B}{\mathcal B}
\newcommand{\CC}{\mathcal C}

\newcommand{\X}{\mathcal X}
\renewcommand{\L}{\mathcal L}

\newcommand{\A}{\mathcal A}

\newcommand{\C}{\mathbb C}

\newcommand{\midtext}[1]{\quad\text{#1}\quad}
\newcommand{\righttext}[1]{\qquad\text{#1 }}

\DeclareMathOperator{\id}{id}

\DeclareMathOperator{\paut}{PAut}

\DeclareMathOperator{\supp}{supp}
\DeclareMathOperator*{\clsp}{\overline{span}}
\DeclareMathOperator*{\spn}{span}

\newcommand{\norm}[1]{\left\| #1 \right\|}

\newcommand{\iso}{\overset{\cong}{\longrightarrow}}

\begin{document}

\title[Actions of groupoids and inverse semigroups]
{$C^*$-actions of $r$-discrete groupoids and inverse semigroups}

\author{John Quigg}
\address{Department of Mathematics\\Arizona State University\\
Tempe, Arizona 85287}
\email[J. Quigg]{quigg@math.la.asu.edu}

\author{N\'andor Sieben}
\email[N. Sieben]{Nandor.Sieben@asu.edu}

\thanks{Research partially supported by National Science Foundation
Grant No. DMS9401253}



\begin{abstract}
Groupoid actions on $C^*$-bundles and inverse 
semigroup actions on $C^*$-algebras are closely related when the 
groupoid is $r$-discrete.
\end{abstract}

\maketitle

\section{Introduction}
Many important $C^{*}$-algebras, such as AF-algebras, 
Cuntz-Krieger algebras, graph algebras and foliation 
$C^{*}$-algebras, are the $C^{*}$-algebras of $r$-discrete groupoids.  
These $C^{*}$-algebras are often associated with inverse 
semigroups through the $C^{*}$-algebra of the inverse 
semigroup \cite{han-rae} or through a crossed 
product construction as in Kumjian's localization \cite{kum:local}.  
Nica \cite{nic:inv} connects groupoid $C^{*}$-algebras with 
the partial 
crossed product $C^{*}$-algebras of Exel \cite{exe:partial} and 
McClanahan \cite{mcc}.  This gives another connection between 
groupoid $C^{*}$-algebras and inverse semigroup $C^{*}$-algebras 
since 
\cite{sie:partial} and \cite{exe:inverse} show that discrete partial
crossed products are basically special cases of the
inverse 
semigroup crossed products of \cite{sie:partial}, \cite{pat:book}, 
\cite{sie:thesis}.  

The heart of these connections is Renault's observation 
in \cite{ren:thesis} that an $r$-discrete groupoid can be 
recovered from the way the inverse semigroup of open 
$G$-sets acts on the unit space of the groupoid. In 
the upcoming \cite{pat:book} Paterson further develops 
this connection by showing that the $C^{*}$-algebra of an 
$r$-discrete groupoid $G$ is the crossed product of $C_0(G^0)$ by 
the action of the inverse semigroup of open $G$-sets.

The purpose of this paper is to explore this connection 
on the level of $C^{*}$-crossed products.  Renault 
\cite{ren:representation} defines a $C^{*}$-action of a groupoid as a 
functor to the category of $C^{*}$-algebras and 
homomorphisms, in which the collection of object 
$C^{*}$-algebras are glued together as a $C^{*}$-bundle over $G^0$ 
and the action is appropriately continuous.  We associate 
to this an action of 
any sufficiently large
inverse semigroup $S$ of open 
$G$-sets on the $C_0$-section algebra of the bundle.  
Conversely, starting with an action (satisfying certain 
mild conditions) of $S$ on a $C^{*}$-algebra $B$, we obtain 
an
associated $C^{*}$-bundle over $G^0$ via the realization that 
$C_0(G^0)$ will act as central multipliers of $B$.  Then we 
construct the groupoid action using the expected `germs 
of local automorphisms'
approach that goes back to \cite{hae} and \cite{rei}.  
The $C^{*}$-bundles arising this way 
are typically only upper semicontinuous, rather than 
continuous.  So we use a slight generalization of 
Renault's theory.  

The philosophy is that inverse 
semigroups and $r$-discrete groupoids are two sides of the 
same coin;  passing back and forth between groupoid and 
inverse semigroup constructions may benefit both 
theories.  The theory of groupoid $C^{*}$-algebras is more 
developed, but the inverse semigroup theory is more algebraic.
For example, 
we can show that the $C^*$-algebra of an $r$-discrete 
groupoid is an enveloping $C^*$-algebra without using Renault's 
disintegration theorem.  
In fact one could suspect that for $r$-discrete groupoids the 
disintegration theorem follows from the less complicated inverse 
semigroup disintegration theorem.  Other applications 
could include inverse semigroup versions of Kumjian's 
\cite{kum:fell} groupoid Fell bundles and Renault's 
imprimitivity theorem \cite{ren:representation} (see also 
\cite{rae:symmetric}).  This latter may be an important 
step towards finding a regular representation for inverse 
semigroup actions, 
which is very much needed for 
coactions and crossed product duality.  

After some preliminary results, we 
introduce a slight generalization of
Renault's groupoid 
actions 
in section 3. In section 4
we recall the basic theory of inverse semigroup actions.
In sections 5 and 6 we show how to pass back and forth
between groupoid and inverse semigroup actions. In section 
7 we prove our main theorem by showing that the crossed products 
of the corresponding groupoid and inverse semigroup actions 
are isomorphic. Finally, as an application, we recover the Hausdorff 
case of
Paterson's theorem connecting groupoid $C^{*}$-algebras and inverse 
semigroup actions. 
Starting with an inverse semigroup, Paterson 
builds a universal groupoid \cite{pat:book}. This groupoid is not Hausdorff in general.
Since we only work with Hausdorff groupoids we did not 
use this universal groupoid, rather we assumed that our inverse
semigroup is always a semigroup of open $G$-sets of a groupoid. It is
likely that our approach works for non-Hausdorff groupoids and the
theory generalizes 
to the level of the universal groupoid. 

\section{Preliminaries}
\label{prelim}
We will need the following elementary results on representations of
$C^*$-algebras. Since we could not find a reference, we include the
proofs for the convenience of the 
reader. When we refer to a `representation' of a
$^*$-algebra, we mean a $^*$-homomorphism of the algebra into the
bounded operators on a Hilbert space.

\begin{defn}
Let $D$ be a $^*$-algebra. We say $D$ \emph{has an enveloping
$C^*$-algebra} if the supremum of the $C^*$-seminorms on $D$ is
finite, and in this case we call the Hausdorff completion of $D$
relative to this largest $C^*$-seminorm the \emph{enveloping
$C^*$-algebra} of $D$.
\end{defn}

\begin{rem}
Thus, if $D$ is a $^*$-subalgebra of a $C^*$-algebra $B$, and if every
representation of $D$ is bounded in the norm inherited from $B$, then
the closure of $D$ in $B$ is the enveloping $C^*$-algebra of $D$.

Conversely, if the closure of $D$ in $B$ is the enveloping $C^*$-algebra
of $D$, then every representation of $D$ is contractive.
\end{rem}

Our first elementary result about enveloping $C^*$-algebras is that
ideals have them.

\begin{lem}
\label{ideal}
Let $I$ be a  two-sided, not necessarily closed, $^*$-ideal of
a $C^*$-subalgebra $B$. Then the closure of $I$ in $B$ is the
enveloping $C^*$-algebra of $I$.
\end{lem}

\begin{proof}
Let $\pi$ be a representation of $I$ on $H$.
Since we can replace $H$
by the closure of $\pi(I)H$, we may as well assume $\pi(I)H$ is dense
in $H$. Then of course $\pi(I^2)H$ is also dense in $H$, so
it suffices to show
\[
\norm{\sum_1^n\pi(ab_ic_i)\xi_i}
\le\norm{a}\norm{\sum_1^n\pi(b_ic_i)\xi_i}
\righttext{for}a,b_1,\dots,b_n,c_1,\dots,c_n\in I,
\xi_1,\dots,\xi_n\in H.
\]
We use the Effros-Hahn trick: put
\[
d=\bigl(\norm{a}^2-a^*a\bigr)^{1/2}\in M(B).
\]
We have
\begin{align*}
\norm{\sum_1^n\pi(ab_ic_i)\xi_i}^2
&=\sum_{i,j}(\pi(a^*ab_ic_i)\xi_i,\pi(b_jc_j)\xi_j)\\
&=\sum_{i,j}(\pi(\norm{a}^2b_ic_i-d^2b_ic_i)\xi_i,\pi(b_jc_j)\xi_j)\\
&=\sum_{i,j}\norm{a}^2(\pi(b_ic_i)\xi_i,\pi(b_jc_j)\xi_j)
-\sum_{i,j}(\pi(d^2b_ic_i)\xi_i,\pi(b_jc_j)\xi_j),\\
\intertext{which makes sense since $M(B)I^2\subset BI\subset I$,}
&=\norm{a}^2\norm{\sum_1^n\pi(b_ic_i)\xi_i}^2
-\norm{\sum_1^n\pi(db_ic_i)\xi_i}^2\\
&\le\norm{a}^2\norm{\sum_1^n\pi(b_ic_i)\xi_i}^2,
\end{align*}
as desired.
\end{proof}

\begin{rem}
As usual with the Effros-Hahn trick, the above argument shows even more:
we only need to assume $\pi$ is a $^*$-homomorphism of $I$ into the
$^*$-algebra of \emph{not necessarily bounded} linear operators on $ H$,
such that
\[
\overline{\pi(I) H}= H
\midtext{and}
(\pi(a)\xi,\eta)=(\xi,\pi(a^*)\eta),
\]
and then the argument shows each $\pi(a)$ is automatically bounded.
\end{rem}

\begin{defn}
Say that a family $\{B_e\}_{e\in E}$ of closed ideals of a $C^*$-algebra
$B$ is \emph{upward-directed} if for all $e,f\in E$ there exists $g\in
E$ such that $B_e\cup B_f\subset B_g$.
\end{defn}

The following elementary result
allows us to paste together consistent representations of 
an upward-directed family of ideals.

\begin{lem}
\label{paste}
Let $\{B_e\}_{e\in E}$ be an upward-directed family of closed ideals
with dense span in 
a $C^*$-algebra $B$, and suppose that for each $e\in E$ we have a
representation $\pi_e$ of $B_e$ on a common Hilbert space $ H$, such that
\begin{enumerate}
\item $\pi_e=\pi_f|{B_e}$ whenever $B_e\subset B_f$, and
\item $\spn_{e\in E}\pi_e(B_e) H$ is dense in $ H$.
\end{enumerate}
Then there is a unique representation of $B$ on $ H$ which extends
every $\pi_e$.
\end{lem}

\begin{proof}
By upward-directedness the union $\bigcup_{e\in E}B_e$ is a dense ideal
of $B$. The consistency condition (i) guarantees that the union of the
$\pi_e$'s is a representation $\pi$ of $\bigcup_{e\in E}B_e$, which is
of course contractive since each $\pi_e$ is. The nondegeneracy condition
(ii) shows $\pi$ extends uniquely to a representation of $B$.
\end{proof}

\begin{rem}
We did not need to appeal to \lemref{ideal} to extend the representation
$\pi$ from $\bigcup_{e\in E}B_e$ to $B$ since the hypotheses already
told us $\pi$ was contractive.
\end{rem}

\section{Groupoid actions}
\label{gpoid}

Throughout, $G$ will be a 
locally compact Hausdorff
groupoid with a Haar system. 
After the next few paragraphs $G$ will be $r$-discrete.
We should comment a little on the still-evolving definition of 
this term.
Renault \cite[Definition 1.2.6]{ren:thesis}
defined a locally compact Hausdorff groupoid $G$ to be $r$-discrete
if the unit space $G^0$ is open in $G$. However, this does not quite
give all that one wants, in particular a Haar system. In fact, Renault
proved that a groupoid which is $r$-discrete in his sense has a Haar
system if and only if the range and domain maps $r$ and $d$ are local
homeomorphisms, and in this case counting measures give a Haar system.
The current
fashion is to build this into the definition of $r$-discrete. Perhaps one
of the most elegant `modern' definitions is Paterson's
\cite{pat:book}: 
for a locally compact (\emph{not necessarily
Hausdorff}) groupoid $G$ Paterson defines $G^{op}$
as the set of open Hausdorff subsets $U$ of $G$ such that $r|U$ and
$d|U$ are homeomorphisms onto open subsets of $G$,
and he calls $G$
$r$-discrete if $G^{op}$ is a base for the topology of $G$.
This is compatible with Renault's definition when $G$ is Hausdorff,
and in this case (which is our primary concern) the $G^{op}$
condition just means the range and domain maps $r$ and $d$ are local
homeomorphisms. Throughout this paper, when we say
$G$ is \emph{$r$-discrete}
we will always mean $G$ is
locally compact and Hausdorff, and
$r$ and $d$ are local homeomorphisms.
In particular, in an $r$-discrete groupoid $G$ the range and domain maps
$r$ and $d$
are open, so the elements of $G^{op}$ are just the open $G$-sets, 
also called \emph{bisections}
(recall that a $G$-set is defined as a subset of $G$ on which $r$ and
$d$ are injective).
It was with some hesitation that we imposed the Hausdorff condition;
many $r$-discrete groupoids occurring in nature are non-Hausdorff, and
it would seem reasonable to expect that our results hold for all such
groupoids. However, our techniques seem to depend upon Hausdorffness;
we intend to explore this in future work.

In \cite{ren:representation} Renault developed a notion of 
actions and crossed products of groupoids in the $C^*$-category, and 
we will 
require 
a slight generalization of his definition of action.
Algebraically, an \emph{action} of $G$ is a functor $\alpha$ from $G$ 
to the category of $C^*$-algebras and homomorphisms, with 
$\alpha_x\colon A_{d(x)}\to A_{r(x)}$ for $x\in G$. Of course, since 
$G$ has inverses, by functoriality each $\alpha_x$ will be an 
isomorphism of $A_{d(x)}$ onto $A_{r(x)}$. Topologically, the 
collection $\A=\{A_u\}_{u\in G^0}$ of $C^*$-algebras must be glued 
together so we can formulate a continuity condition. Renault 
requires $\A$ to be a continuous $C^*$-bundle, 
but we will relax this 
to \emph{upper} semicontinuity. So, we have a continuous, open 
surjection of $\A$ onto $G^0$, and the norm function $\norm{\cdot}$ on 
$\A$ is only upper semicontinuous (see \cite{bla}, \cite{dix}, 
\cite{dup-gil}, 
\cite{nil:bundle}, 
and \cite{rie:bundle}).
Pulling back via the domain map $d\colon G\to G^0$, we
get a Banach bundle $d^*\A=\{(x,a)\in G\times \A:a\in A_{d(x)}\}$ over
$G$.  The continuity condition for $\alpha$ is that the map
$(x,a)\mapsto\alpha_x(a)$ from $d^*\A$ to $\A$ be continuous.
For his main results, Renault also requires separability assumptions,
partly because he appeals to direct integral theory. Since we will not need
direct integrals, we can dispense here with separability hypotheses.

We often need to work with sections, rather than the bundle itself. 
Recall that upper semicontinuity of the $C^*$-bundle $\A$ can be 
equivalently expressed as the condition that for every
$f$ in the $C_0$-section algebra $\Gamma_0(\A)$ 
the 
map $u\mapsto\norm{f(u)}$ is upper semicontinuous. We will need the 
following characterization of continuity for actions in terms of 
compactly supported
sections.

\begin{lem}
\label{continuous}
Let $\A$ be an upper semicontinuous $C^*$-bundle over $G^0$, and 
let $\alpha$ be a functor from $G$ to the $C^*$-category with 
$\alpha_x\colon A_{d(x)}\to A_{r(x)}$ for each $x\in G$. Then 
$\alpha$ is continuous \textup(hence an action\textup) if and only 
if for all $x\in G$, $a\in A_{d(x)}$, and $f,g\in\Gamma_c(\A)$ with 
$f(r(x))=\alpha_x(a)$ and $g(d(x))=a$,
\[
\lim_{y\to x}\norm{f(r(y))-\alpha_y(g(d(y)))}=0.
\]
\end{lem}
\begin{proof}
First assume $\alpha$ is continuous, and fix $x,a,f,g$ as above. As 
$y\to x$ we have $r(y)\to r(x)$, so $f(r(y))\to 
f(r(x))=\alpha_x(a)$ by continuity. Similarly, $g(d(y))\to 
g(d(x))=a$, so $\alpha_y(g(d(y)))\to\alpha_x(a)$ by continuity of 
$\alpha$. Hence
\[
0=\norm{f(r(x))-\alpha_x(a)}
\ge\limsup_{y\to x}\norm{f(r(y))-\alpha_y(g(d(y)))}
\]
by upper semicontinuity, giving 
$\norm{f(r(y))-\alpha_y(g(d(y)))}\to 0$.

Conversely, let $(x,a)\in d^*\A$, and let $V$ be a neighborhood of 
$\alpha_x(a)$ in $\A$. Pick $f,g\in\Gamma_c(\A)$ with 
$f(r(x))=\alpha_x(a)$ and $g(d(x))=a$. By \cite[page 10]{dup-gil} we 
may assume 
\[
V=\bigcup_{u\in U}\{b\in A_u:
\norm{f(u)-b}<\epsilon\}
\]
for some neighborhood $U$ of $r(x)$ in $G^0$ and some 
$\epsilon>0$. Assuming
\[
\lim_{y\to x}\norm{f(r(y))-\alpha_y(g(d(y)))}=0,
\]
we can find a neighborhood $N$ of $x$ in $G$ such that $r(N)\subset 
U$ and
\[
\norm{f(r(y))-\alpha_y(g(d(y)))}<\frac{\epsilon}{2}
\righttext{for}y\in N.
\]
Put
\[
W=\bigcup_{u\in d(N)}\{b\in A_u:
\norm{g(u)-b}<\frac{\epsilon}{2}\}.
\]
Then $N*W:=(N\times W)\cap d^*\A$
is a neighborhood of $(x,a)$ in $d^*\A$, and for $(y,b)\in 
N*W$ we have
\begin{align*}
&\norm{f(r(y))-\alpha_y(b)}\\
&\quad\le\norm{f(r(y))-\alpha_y(g(d(y)))}
+\norm{\alpha_y(g(d(y)))-\alpha_y(b)}\\
&\quad<\frac{\epsilon}{2}+\norm{g(d(y))-b}
<\epsilon,
\end{align*}
so $\alpha_y(b)\in V$, and we have shown $(y,b)\mapsto\alpha_y(b)$ is 
continuous at $(x,a)$.
\end{proof}

Now let $\alpha$ be an action of $G$ on $\A$, and let
$r^*\A=\{(a,x)\in\A\times G:a\in A_{r(x)}\}$ be the pull-back via the
range map $r\colon G\to G^0$.  Renault \cite{ren:representation} shows that the
vector space $\Gamma_c(r^*\A)$ of compactly supported continuous sections
becomes a $^*$-algebra with the operations
\begin{gather*}
(fg)(x)=\int f(y)\alpha_y(g(y^{-1}x))\,d\lambda^{r(x)}(y)
\midtext{and}
f^*(x)=\alpha_x(f(x^{-1}))^*.
\end{gather*}
He then defines 
(in 
\cite[beginning of Section 5]{ren:representation})
the 
\emph{crossed product} $\A\times_\alpha G$ as the 
Hausdorff completion 
of $\Gamma_c(r^*\A)$ in the supremum of the $C^*$-seminorms
\[
f\mapsto\norm{\Pi(f)},
\]
where $\Pi$ runs over all representations of $\Gamma_c(r^*\A)$ which
are continuous from the inductive limit topology to the weak operator
topology.
He shows as a consequence of his decomposition theorem \cite[Th\'eor\`eme
4.1]{ren:representation} that this supremum is finite. We will give
an independent proof of this for \emph{$r$-discrete} groupoids,
as an application of our isomorphism (see \thmref{bundle}) between
$\Gamma_c(r^*\A)$ and a $^*$-algebra associated with a corresponding
inverse semigroup action.  In fact, this will show that for $r$-discrete
groupoids the crossed product is the \emph{enveloping $C^*$-algebra} of
$\Gamma_c(r^*\A)$, because the inductive limit continuity is automatic,
as we show in the following proposition.

\begin{prop}
\label{G envelope}
If $\alpha$ is an action of an $r$-discrete groupoid $G$ 
on an upper semicontinuous $C^*$-bundle $\A$, then
every representation
of $\Gamma_c(r^*\A)$ is continuous from the inductive limit topology to
the weak operator topology.
\end{prop}

\begin{proof}
Let $\Pi$ be a representation of $\Gamma_c(r^*\A)$ on a Hilbert space
$ H$.
For each subset $T$ of $G$ define
\[
\Gamma_T(r^*\A)=\{f\in\Gamma_c(r^*\A):\supp f\subset T\}.
\]

Claim: it suffices to show that for each $s\in G^{op}$
the restriction $\Pi|_{\Gamma_s(r^*\A)}$ is continuous from the (sup)
norm topology to the weak operator topology. To see this, let $K$ be a
compact subset of $G$, and take $\xi,\eta\in H$. Find $s_1,\dots,s_n\in
G^{op}$ such that $K\subset\bigcup_1^n s_i$, and choose a partition of unity
$\{\phi_i\}_1^n$ subordinate to the open cover $\{s_i\}_1^n$ of $K$,
so that $\supp\phi_i\subset s_i$ and $\sum_1^n\phi_i=1$ on $K$.
For
$f\in\Gamma_K(r^*\A)$ we have
\[
(\Pi(f)\xi,\eta)
=\biggl(\Pi\Bigl(\sum_1^n f\phi_i\Bigr)\xi,\eta\biggr)
=\sum_1^n(\Pi(f\phi_i)\xi,\eta).
\]
Now just observe that the map $f\mapsto
f\phi_i\colon\Gamma_K(r^*\A)\to\Gamma_{s_i}(r^*\A)$ is norm continuous, and
the claim follows.

So, fix $s\in G^{op}$. We will in fact show $\Pi|_{\Gamma_s(r^*\A)}$ is continuous
for the \emph{norm} topologies of $\Gamma_s(r^*\A)$ and $\L( H)$. Note that
for $f\in\Gamma_s(r^*\A)$ we have
\[
f^*f(x)=\sum_{r(y)=r(x)}\alpha_y\bigl(f(y^{-1})^*f(y^{-1}x)\bigr),
\]
and for any nonzero term in this sum we have $y\in s^*$
and $x\in ys\subset s^*s=d(s)$. Hence, the product $f^*f$
is in $\Gamma_{d(s)}(r^*\A)$, and for each $x\in d(s)$ we have
$f^*f(x)=\alpha_y(f(y)^*f(y))$, 
where $y$ is the unique element of $s$ with
$d(y)=x$.  
Now, $\Gamma_{d(s)}(r^*\A)$ is an ideal in the $C_0$-section
algebra $\Gamma_0((r^*\A)|_{G^0})$, which in turn is a $C^*$-subalgebra
of $\A\times_\alpha G$.  Consequently, \lemref{ideal} tells us the
restriction $\Pi|_{\Gamma_{d(s)}(r^*\A)}$ is contractive, so
\[
\norm{\Pi(f)}^2
\le\norm{f^*f}_\infty
=\norm{f}_\infty^2,
\]
since
\[
\sup_{x\in d(s)}\norm{f^*f(x)}
=\sup_{y\in s}\norm{f(y)^*f(y)}
=\sup_{y\in s}\norm{f(y)}^2.
\]
\end{proof}

\section{Inverse semigroup actions}

Here we mainly follow the conventions of 
\cite{pat:book}, \cite{sie:partial}, and
\cite{sie:thesis}.
Let $B$ be a $C^*$-algebra. A \emph{partial
automorphism} \cite{exe:partial} of $B$ is an isomorphism between
two (closed) ideals of $B$. The set of all partial automorphisms
of $B$ forms an inverse semigroup $\paut B$ under composition. An
\emph{action} of an inverse semigroup $S$ on $B$ is a homomorphism
$\beta\colon S\to\paut B$ such that the domain ideals of the partial
automorphisms $\{\beta_s\}_{s\in S}$ 
are upward-directed and have dense span in $B$.
Of course, for $s\in S$ the domain of the partial automorphism
$\beta_s$ will also be the range of $\beta_{s^*}$, and furthermore
will only depend upon the domain idempotent $d(s):=s^*s$. Thus, for
each \emph{idempotent} 
$e\in E_S$ we have an ideal $B_e$ of $B$, and each $\beta_s$ is an
isomorphism of $B_{d(s)}$ onto $B_{r(s)}$.  
Recall from \secref{prelim} that the upward-directedness of
the ideals $\{B_e\}_{e\in E_S}$ means that any two of them are contained
in a third,
and is automatic
if
$E_S$ itself is 
upward-directed in the sense that for all $e,f\in E_S$ there exists
$g\in E_S$ such that $e,f\le g$.

A \emph{representation} of $S$ on a Hilbert space $H$ is a 
$^*$-homomorphism $U$ of $S$ to $\L(H)$ such that the span of the 
ranges of the operators $\{U_s\}_{s\in S}$ is dense in $H$. Of 
course, each $U_s$ will be a partial isometry, since
\[
U_sU_s^*U_s=U_sU_{s^*}U_s=U_{ss^*s}=U_s.
\]

A \emph{covariant representation} of an inverse semigroup action 
$(B,S,\beta)$ is a pair $(\pi,U)$ consisting of a nondegenerate 
representation $\pi$ of $B$ and a representation $U$ of $S$, both on 
the same Hilbert space $H$, such that 
\[
U_eH=\pi(B_e)H\righttext{for}e\in E_S
\]
and
\[
U_s\pi(b)U_s^*=\pi(\beta_s(b))\righttext{for}b\in B_{d(s)}.
\]
Note that when we are checking whether a pair $(\pi,U)$ is
covariant, we do not need to verify the span of the ranges of the $U_s$
is dense, since this follows from nondegeneracy of $\pi$ and the
covariance condition.

The disjoint union of $\{B_e\}_{e\in E_S}$ forms a Banach bundle $\B$ over
the discrete space $E_S$. Pulling back via the range map $r\colon S\to
E_S$, we get a Banach bundle $r^*\B=\{(b,s)\in\B\times S:b\in B_{r(s)}\}$
over $S$.  The set $\Gamma_c(r^*\B)$ of finitely supported sections
becomes a $^*$-algebra with operations 
defined on the generators by
\[
(b,s)(c,t)=(\beta_s(\beta_{s^*}(b)c),st)\midtext{and}
(b,s)^*=(\beta_{s^*}(b^*),s^*),
\]
and then extended additively.  For every covariant representation
$(\pi,U)$ of $(B,S,\beta)$ the \emph{integrated form} of $(\pi,U)$
is the representation $\Pi$ of $\Gamma_c(r^*\B)$ defined by
\[
\Pi\biggl(\sum_1^n(b_i,s_i)\biggr)=\sum_1^n\pi(b_i)U_{s_i}.
\]
The integrated form $\Pi$ is nondegenerate, and we have
\[
\Pi(\Gamma_c(r^*\B))=\spn_{s\in S}\pi(B_{r(s)})U_s
\midtext{and}
\norm{\Pi\biggl(\sum_1^n(b_i,s_i)\biggr)}\le\sum_1^n\norm{b_i}.
\]
The \emph{crossed product} of $(B,S,\beta)$ is the Hausdorff completion 
$B\times_\beta S$ of $\Gamma_c(r^*\B)$ in the supremum of the 
$C^*$-seminorms
\[
f\mapsto\norm{\Pi(f)},
\]
where $\Pi$ runs over the integrated forms of all covariant
representations.  Warning: there is some collapsing when $\Gamma_c(r^*\B)$
is mapped into $B\times_\beta S$, since whenever $\Pi$ is the integrated
form of a covariant representation we have
\[
\Pi(b,s)=\Pi(b,t)\righttext{for}b\in B_{r(s)},s\le t
\]
(see
\cite[Proposition 3.3.2]{pat:book}, \cite[Lemma 3.4.4]{sie:thesis},
\cite[Lemma 4.5]{sie:partial}).
For $b\in B_{r(s)}$ let $[b,s]$
denote the image of $(b,s)$ in $B\times_\beta S$, so that
\[
B\times_\beta S=\clsp\{[b,s]:b\in B_{r(s)},s\in S\}.
\]
For every covariant representation $(\pi,U)$ there is a unique
representation $\pi\times U$ of $B\times_\beta S$ such that
\[
(\pi\times U)[b,s]=\pi(b)U_s\righttext{for}b\in B_{r(s)},s\in S,
\]
and in fact this gives a bijection between the covariant representations
of the action $(B,S,\beta)$ and the nondegenerate representations of
the $C^*$-algebra $B\times_\beta S$ 
\cite[Corollary 3.3.1]{pat:book}, \cite[Proposition 3.4.7]{sie:thesis}.

Caution: the crossed product $B\times_\beta S$ is (usually) \emph{not} the
enveloping $C^*$-algebra of $\Gamma_c(r^*\B)$; although it is true (and not
hard to show) that $\Gamma_c(r^*\B)$ does in fact have an enveloping
$C^*$-algebra, there can in general be representations of $\Gamma_c(r^*\B)$
which are not integrated forms of covariant representations \cite[Example
4.9]{sie:partial}. Which representations of $\Gamma_c(r^*\B)$ are integrated
forms?  Paterson 
\cite{pat:book} 
has found an answer:

\begin{defn}
A representation $\Pi$ of $\Gamma_c(r^*\B)$ is called \emph{coherent} if
\[
\Pi(b,e)=\Pi(b,f)\righttext{whenever}b\in B_e \text{ and } e\le f.
\]
\end{defn}

Actually, Paterson 
requires $\Pi(b,e)=\Pi(b,f)$ whenever $b\in B_eB_f$, which is clearly
equivalent to the above logically weaker condition, and he builds this
condition right into his definition
of a representation of $\Gamma_c(r^*\B)$. 
He also requires the map $\Pi$ to be bounded 
in the $L^1$-norm
\[
\norm{\sum(b_s,s)}_1:=\sum\norm{b_s}\left(=\sum\norm{(b_s,s)}\right),
\]
although 
this follows automatically from \propref{coherent} below. Clearly,
$\Pi$ is coherent if and
only if its kernel contains the
ideal generated by the subspace
\[
\spn\{(b,e)-(b,f):b\in B_e,e,f\in E_S,e\le f\}.
\]
We need to 
identify this ideal explicitly. Recall that the partial order in an
inverse semigroup is given by
\[
s\le t\midtext{if and only if}s=ss^*t,
\]
which should be regarded as saying $s$ is a `restriction' of $t$.

\begin{lem}
The ideal of $\Gamma_c(r^*\B)$ generated by
$\{(b,e)-(b,f):b\in B_e,e,f\in E_S,e\le f\}$ 
coincides with the subspace
\[
I_\beta:=\spn\{(b,s)-(b,t):b\in B_s,s\le t\}.
\]

\end{lem}

\begin{proof}
It is easy to check that $I_\beta$ is a self-adjoint left ideal. Hence,
it suffices to observe that for $b\in B_{r(s)}$ and $s\le t$ we can
factor $b=cd$ for some $c,d\in B_{r(s)}$, and then
\begin{align*}
(b,s)-(b,t)
&=(cd,ss^*t)-(cd,tt^*t)
=(c,ss^*)(d,t)-(c,tt^*)(d,t)
\\&=\bigl((c,ss^*)-(c,tt^*)\bigr)(d,t).
\end{align*}
\end{proof}

Is there a coherent representation whose kernel is $I_\beta$? We do not
know in general, but it follows from \thmref{bundle} below that the
answer is yes for representations related to groupoids.

The following 
proposition, which appears in various forms in 
\cite[Corollary 3.3.1]{pat:book},
and \cite[Proposition 3.4.7]{sie:thesis},
establishes a bijective correspondence between coherent, 
nondegenerate representations of $\Gamma_c(r^*\B)$ and covariant
representations of $(B,S,\beta)$. 
We include the outline of the argument for the convenience of the
reader; in particular, this is the first time the automatic continuity
of representations of $\Gamma_c(r^*\B)$ has been adequately handled.

\begin{prop}[\cite{pat:book}, \cite{sie:thesis}]
\label{coherent}
Every representation of $\Gamma_c(r^*\B)$ is contractive on each fiber
$(B_{r(s)},s)$.
A nondegenerate representation $\Pi$ of $\Gamma_c(r^*\B)$ is the
integrated form of a covariant representation of $(B,S,\beta)$ if and
only if $\Pi$ is coherent.
\end{prop}

\begin{proof}
Let $\Pi$ be a representation of $\Gamma_c(r^*\B)$.
For each $e\in E_S$ define a
representation $\pi_e$ of $B_e$ by
\[
\pi_e(b)=\Pi(b,e).
\]
The first statement follows from the estimate
\begin{align*}
\norm{\Pi(b,s)}^2
&=\norm{\Pi(b,s)^*\Pi(b,s)}
=\norm{\Pi(\beta_{s^*}(b^*),s^*)\Pi(b,s)}
\\&=\norm{\Pi\bigl((\beta_{s^*}(b^*),s^*)(b,s)\bigr)}
=\norm{\Pi\bigl(\beta_{s^*}(\beta_s\beta_{s^*}(b^*)b),s^*s\bigr)}
\\&=\norm{\Pi(\beta_{s^*}(b^*b),d(s))}
=\norm{\pi_{d(s)}(\beta_{s^*}(b^*b))}
\\&\le\norm{\beta_{s^*}(b^*b)}
=\norm{b^*b}
=\norm{b}^2
=\norm{(b,s)}^2.
\end{align*}
For the other part, first
assume $\Pi$ is the integrated form of a covariant representation
$(\pi,U)$. 
Then for $b\in B_{e}$ and $e\le f$ in $E_S$,
\begin{align*}
\Pi(b,e)
&=\pi(b)U_e
=\pi(b)U_{ef}
=\pi(b)U_eU_f
=\pi(b)U_f
=\Pi(b,f),
\end{align*}
since $U_e H=\pi(B_{e}) H$, and so $\Pi$ is coherent.

Conversely, assume $\Pi$ is 
a nondegenerate and coherent representation of $\Gamma_c(r^*\B)$ on
$ H$, and recall the above definition of the $\pi_e$. We have
$\pi_e=\pi_f|{B_e}$ whenever $B_e\subset B_f$, 
by coherence,
and
$\spn_e\pi_e(B_e) H$ is dense in $ H$, so by \lemref{paste} there is a
unique representation $\pi$ of $B$ on $ H$ which extends every $\pi_e$.

To get the other half of our covariant representation we use the
construction of McClanahan \cite[Proposition 2.8]{mcc} (where the context
was partial actions of groups), which was adapted to inverse semigroup
actions by the second author \cite[proof of Proposition 4.7]{sie:partial}.
Fix $s\in S$ and a bounded approximate identity $\{b_i\}$ for
$B_{r(s)}$. Claim: the net $\{\Pi(b_i,s)\}$ is Cauchy in the strong
operator topology. To see this, take $\xi\in H$. Then
\begin{align*}
\norm{\bigl(\Pi(b_i,s)-\Pi(b_j,s)\bigr)\xi}^2
&=\Bigl(\bigl(\Pi(b_i,s)-\Pi(b_j,s)\bigr)^*
\bigl(\Pi(b_i,s)-\Pi(b_j,s)\bigr)\xi,\xi\Bigr)
\\&=\biggl(\Pi\Bigl(\bigl((b_i,s)-(b_j,s)\bigr)^*
\bigl((b_i,s)-(b_j,s)\bigr)\Bigr)\xi,\xi\biggr)
\\&=\Bigl(\Pi\bigl((\beta_{s^*}(b_i),s^*)(b_i,s)
+(\beta_{s^*}(b_j),s^*)(b_j,s)
\\&\phantom{\Bigl(\Pi\bigl(}\quad-(\beta_{s^*}(b_i),s^*)(b_j,s)
-(\beta_{s^*}(b_j),s^*)(b_i,s)\bigr)\xi,\xi\Bigr)
\\&=\Bigl(\pi_{d(s)}
\bigl(\beta_{s^*}(b_i^2+b_j^2-b_ib_j-b_jb_i)\bigr)\xi,\xi\Bigr)
\\&\to 0,
\end{align*}
since $\{b_i^2\}$, $\{b_j^2\}$, $\{b_ib_j\}$, and $\{b_jb_i\}$ are all
bounded approximate identities (the latter two for the product
direction) for $B_{r(s)}$, hence their images under $\beta_{s^*}$ are
approximate identities for $B_{d(s)}$. This verifies the claim, so we
can let $U_s$ be the strong operator limit of $\{\Pi(b_i,s)\}$.

Since multiplication and involution are jointly strong operator
continuous on bounded sets,
\begin{align}
\label{isom}
\begin{split}
U_s^*U_s
&=\lim\Pi(b_i,s)^*\Pi(b_i,s)
=\lim\Pi(\beta_{s^*}(b_i^2),d(s))
\\&=\lim\pi_{d(s)}(\beta_{s^*}(b_i^2))
\\&=\text{projection on }\pi_{d(s)}(B_{d(s)}) H.
\end{split}
\end{align}
Hence, $U_s$ is a partial isometry with initial subspace
$\pi_{d(s)}(B_{d(s)}) H$. In particular, the following computation shows
$U_s$ is independent of the choice of the bounded approximate identity
$\{b_i\}$: for $c\in B_{d(s)}$, $\xi\in H$
\begin{align}
\label{indep}
\begin{split}
U_s\pi_{d(s)}(c)\xi
&=\lim\Pi(b_i,s)\Pi(c,d(s))\xi
=\lim\Pi(\beta_s(\beta_{s^*}(b_i)c),s)\xi
\\&=\Pi(\beta_s(c),s)\xi,
\end{split}
\end{align}
since $\{\beta_{s^*}(b_i)\}$ is a bounded approximate identity for
$B_{d(s)}$ and $\Pi$ is continuous on the fiber $(B_{r(s)},s)$ of
$r^*\B$.

To see that $U$ is multiplicative, take bounded approximate identities
$\{b_i\}$ for $B_{r(s)}$ and $\{c_j\}$ for $B_{r(t)}$:
\begin{align*}
U_sU_t
&=\lim\Pi(b_i,s)\Pi(c_j,t)
=\lim\Pi(\beta_s(\beta_{s^*}(b_i)c_j),st)
=U_{st},
\end{align*}
since $\{\beta_s(\beta_{s^*}(b_i)c_j)\}$ is a bounded approximate identity
for $\beta_s(B_{d(s)}B_{r(t)})=B_{r(st)}$.  A similar computation shows
$U_s^*=U_{s^*}$.

For covariance, the computation \eqref{isom} implies $U_e H=\pi_e(B_e) H$
for all $e\in E_S$, and for $b\in B_{d(s)}$ the computation \eqref{indep}
shows
\begin{align}
\label{cov}
\begin{split}
U_s\pi(b)
&=\Pi(\beta_s(b),s)
=\Pi(\beta_{s^*}(\beta_s(b^*)),s^*)^*
\\&=\bigl(U_{s^*}\pi(\beta_s(b^*))\bigr)^*
=\pi(\beta_s(b^*))^*U_{s^*}^*
=\pi(\beta_s(b))U_s.
\end{split}
\end{align}

$\pi$ is nondegenerate, since if $\pi(B)\xi=0$ then for all $s\in
S$, $b\in B_{d(s)}$ we have
\[
0=U_s\pi_{d(s)}(b)\xi=\Pi(\beta_s(b),s)\xi,
\]
so $\xi=0$ by nondegeneracy of $\Pi$. Finally, 
the computation
\eqref{cov} also implies 
$\Pi(b,s)=\pi(b)U_s$ for $b\in B_{r(s)}$, so
$\Pi$ is the integrated form of $(\pi,U)$.
\end{proof}

The above proposition allows us to express the crossed product as an
enveloping $C^*$-algebra:

\begin{cor}
The crossed product $B\times_\beta S$ is the enveloping $C^*$-algebra of
the $^*$-algebra $\Gamma_c(r^*\B)/I_\beta$.
\end{cor}

We will actually need a technical generalization of the above
result. Suppose that for each $e\in E_S$ we have a dense ideal $B_e'$
of $B_e$, such that
\[
\beta_s(B_{d(s)}')=B_{r(s)}'
\righttext{for all}s\in S.
\]
Write $\Gamma_c(r^*\B')$ for the linear span of $\{(b,s)\in\B\times
S:b\in B_{r(s)}'\}$ in $\Gamma_c(r^*\B)$. Note that $\Gamma_c(r^*\B')$
is a $^*$-subalgebra (in fact, an ideal) of $\Gamma_c(r^*\B)$ which is
dense for the pointwise convergence topology. We need to know that the
$^*$-algebra $\Gamma_c(r^*\B')$ determines the $C^*$-algebra
$B\times_\beta S$.

\begin{lem}
Every representation $\Pi$ of $\Gamma_c(r^*\B')$ is continuous from the
pointwise convergence
topology to the norm topology of operators, and hence has a
unique extension to a representation $\overline\Pi$ of $\Gamma_c(r^*\B)$.
Moreover, $\Pi$ is coherent if and only if $\overline\Pi$ is.
\end{lem}

\begin{proof}
For the first part it
suffices to show $\Pi$ is norm continuous on each fiber $(B_{r(s)}',s)$
of $r^*\B'$.
Each fiber $(B_{r(s)},s)$ of $r^*\B$ is given the Banach space structure
of $B_{r(s)}$.
For every $e\in E_S$ define a representation $\pi_e$ of
$B_e'$ by $\pi_e(b)=\Pi(b,e)$ for $b\in B_e'$. For $b\in B_{r(s)}'$ we
have
\begin{align*}
\norm{\Pi(b,s)}^2
&=\norm{\Pi(b,s)^*\Pi(b,s)}
=\norm{\Pi\bigl((b,s)^*(b,s)\bigr)}
\\&=\norm{\Pi\bigl((\beta_{s^*}(b^*),s^*)(b,s)\bigr)}
=\norm{\Pi(\beta_{s^*}(b^*b),s^*s)}
\\&=\norm{\pi_{s^*s}(\beta_{s^*}(b^*b))}
\le\norm{\beta_{s^*}(b^*b)},
\\\intertext{by \lemref{ideal}, since
$B_{d(s)}'$ is an ideal of $B_{d(s)}$. So}
\|\Pi(b,s)\|^2&\leq \|\beta_{s^*}(b^*b)\|=\norm{b^*b}
=\norm{b}^2.
\end{align*}

For the coherence, obviously $\Pi$ is coherent if $\overline\Pi$ is, so
assume $\Pi$ is coherent. Take $b\in B_{r(s)}$ and $s\le t$, and choose
a bounded approximate identity $\{c_i\}$ for $B_{r(s)}$ which is
contained in the dense ideal $B_{r(s)}'$. Then
\[
\overline\Pi(b,s)
=\lim_i\Pi(c_ib,s)
=\lim_i\Pi(c_ib,t)
=\overline\Pi(b,t).
\]
\end{proof}

\begin{cor}
\label{S envelope}
With the above notation, $B\times_\beta S$ is the enveloping
$C^*$-algebra of the quotient
\[
\Gamma_c(r^*\B')\bigm/\bigl(I_\beta\cap\Gamma_c(r^*\B')\bigr).
\]
\end{cor}

\begin{proof}
This follows from the above lemma and \propref{coherent}, since a
representation of $\Gamma_c(r^*\B')$ is coherent if and only if it kills
$I_\beta\cap\Gamma_c(r^*\B')$.
\end{proof}

\section{From groupoids to inverse semigroups}

Let $\alpha$ be an action of an $r$-discrete groupoid $G$ on an upper
semicontinuous $C^*$-bundle $\A$.  
Recall that since 
$G$ is $r$-discrete, the family
$G^{op}$ of open $G$-sets is a base for the topology of $G$. Further, 
$G^{op}$
is an inverse semigroup with operations
\[
st=\{xy:(x,y)\in (s\times t)\cap G^2\}\midtext{and}
s^*=\{x^{-1}:x\in s\}.
\]
Note that an element $s$ of 
$G^{op}$ 
has domain idempotent
\begin{align*}
d(s)
&=s^*s
=\{x^{-1}y:(x^{-1},y)\in (s^{-1}\times s)\cap G^2\}\\
&=\{x^{-1}y:x,y\in s,r(x)=r(y)\}\\
&=\{x^{-1}x:x\in s\}
=\{d(x):x\in s\},
\end{align*}
and similarly $s$ has range idempotent
\[
r(s)=ss^*=\{r(x):x\in s\}.
\]
We want to associate to the groupoid action $(\A,G,\alpha)$ an inverse
semigroup action $(B,S,\beta)$. 
For $B$ we take $\Gamma_0(\A)$.
To construct partial automorphisms of $B$, we will
need the following elementary lemma, which is an
easy consequence of, for example, \cite[Proposition 1.4]{dup-gil}.

\begin{lem}
Let $C$ and $D$ be upper semicontinuous $C^*$-bundles over 
locally compact Hausdorff spaces $X$ and $Y$, respectively. 
Let $\phi$ be a homeomorphism of $X$ onto $Y$, and for each 
$x\in X$ let $\gamma_x$ be an isomorphism of $C_x$ onto 
$D_{\phi(x)}$. For $f\in\Gamma_c(C)$ and $y\in Y$ define
\[
\gamma(f)(y)=\gamma_{\phi^{-1}(y)}(f(\phi^{-1}(y)))
\in D_y.
\]
If $\gamma(\Gamma_c(C))\subset\Gamma_0(D)$, then $\gamma$ extends 
uniquely to an isomorphism of $\Gamma_0(C)$ onto $\Gamma_0(D)$. 
Moreover, the extension is given by the above formula for 
$f\in\Gamma_0(C)$.
\end{lem}

For $S$ we want to allow some flexibility; roughly speaking, we can take
any sufficiently large inverse subsemigroup of $G^{op}$.

\begin{defn}
\label{full}
We call an inverse subsemigroup $S$ of $G^{op}$ \emph{full} if $S$ is a
base for the topology of $G$ and $E_S$ is upward-directed in the sense
that every two elements of $E_S$ have a common upper bound.
\end{defn}

Take $S$ to be any inverse subsemigroup of $G^{op}$ which is full in the
above sense. It might be useful to mention that to determine whether
$S$ is a base it is enough to check the idempotents; more precisely,
an inverse subsemigroup $S$ of $G^{op}$ is a base for the topology of
$G$ if and only if $S$ covers $G$ and the idempotent semilattice $E_S$
is a base for the topology of the unit space $G^0$.

Note that $G^{op}$ is a full inverse subsemigroup of itself since it is a 
base for the topology of $G$, and the open $G$-sets in $G^0$ are just
the open sets in $G^0$.  
For a more interesting example, consider a transformation groupoid $G=X\times 
H$ where $H$ is a discrete group with identity $e$. If $\B$ is an 
upward-directed base for the topology of the locally compact space $X$ then 
$S=\{U\times \{h\}:U\in {\B}, h\in H \}$ is full since 
it covers $G$ and $E_S=\{U\times\{e\}:U\in {\B} \}$ is an upward-directed 
base for the topology of $G^0=X\times \{e\}$. 


The ideals of the semigroup action $\beta$ will be given by
\[
B_e=\{f\in\Gamma_0(\A):f(x)=0 \text{ for } x\not\in e\}
\righttext{for}e\in E_S.
\]
The reader can immediately verify that each 
$B_e$ is a closed ideal of $\Gamma_0(\A)$, the span of these
ideals is dense in $\Gamma_0(\A)$, and the family $\{B_e\}_{e\in E_S}$
is 
upward-directed.

\begin{thm}
\label{G-S}
Let $\alpha$ be an action of an $r$-discrete 
groupoid $G$ on an upper
semicontinuous $C^*$-bundle
$\A$, and let $S$ be 
a full inverse semigroup of open $G$-sets. Then there
is a unique action $\beta$ of $S$ on $\Gamma_0(\A)$ such that
\[
\beta_s(f)(u)
=\begin{cases}\alpha_{us}(f(s^*us))&\text{if }u\in r(s),\\
0&\text{else},\end{cases}
\]
for $s\in S$, $f\in B_{d(s)}$, and $u\in G^0$.
\end{thm}

\begin{proof}
We first show that the above formula defines an isomorphism $\beta_s\colon
B_{d(s)}\iso B_{r(s)}$, equivalently, an isomorphism of the $C_0$-section
algebra $\Gamma_0(\A|_{d(s)})$ of the restricted bundle $\A|_{d(s)}$
onto $\Gamma_0(\A|_{r(s)})$, and for this we aim to apply the above
lemma. The map $u\mapsto sus^*$ gives a homeomorphism of $d(s)$ onto
$r(s)$, with inverse $u\mapsto s^*us$. Moreover, $u\mapsto us=s(s^*us)$
is a homeomorphism of $r(s)$ onto $s$, and $d(us)=s^*us$, and similarly
for $u\mapsto su\colon d(s)\to s$. For each $u\in d(s)$, $\alpha_{su}$ is
an isomorphism of $A_u$ onto $A_{sus^*}$. Thus, the proposition follows
from the above lemma once we verify that if $f\in\Gamma_c(\A|_{d(s)})$
then $\beta_s(f)\in\Gamma_0(\A|_{r(s)})$. The continuity properties of
$f$ and $\alpha$ ensure that $\beta_s(f)$ is a continuous section. Also,
if $u\notin s(\supp f)s^*$ then $\beta_s(f)(u)=0$, so $\beta_s(f)$
has compact support.

It remains to check
that $\beta$ is a homomorphism. For
$s,t\in S$ the domain of $\beta_s\beta_t$ is
\begin{align*}
\beta_t^{-1}\bigl( B_{d(s)}\cap B_{r(t)}\bigr)
=\beta_{t^*}\bigl( B_{d(s)\cap r(t)}\bigr)
= B_{t^*d(s)r(t)t}
= B_{t^*s^*st}
= B_{d(st)},
\end{align*}
which is the domain of $\beta_{st}$. For $f\in B_{d(st)}$ 
and $u\in r(st)$ we have $u\in r(s)$ and $s^*us\in r(t)$, so
\begin{align*}
\beta_s\beta_t(f)(u)
&=\alpha_{us}\bigl(\beta_t(f)(s^*us)\bigr)
=\alpha_{us}\bigl(\alpha_{s^*ust}(f(t^*s^*ust))\bigr)
\\&=\alpha_{uss^*ust}(f(t^*s^*ust))
=\beta_{st}(f)(u),
\end{align*}
since
\[
uss^*ust=uust=ust,
\]
and this is enough to show $\beta_s\beta_t=\beta_{st}$.
\end{proof}

\section{From inverse semigroups to groupoids}

As in the preceding section, let $G$ be an $r$-discrete groupoid, and
let $S$ be a full inverse semigroup of open $G$-sets.
Suppose we are given
an action $\beta$ of $S$ on a $C^*$-algebra $B$.  We want to construct
an action of $G$ from which $(B,S,\beta)$ arises as in the construction
of the preceding section. We first need to find an upper semicontinuous
$C^*$-bundle $\A$ over $G^0$ such that $B\cong\Gamma_0(\A)$. We know
from \cite{bla}, \cite{dup-gil}, 
\cite{nil:bundle}, 
and \cite{rie:bundle}
(for example) that this is equivalent to $B$ being a $C_0(G^0)$-algebra,
that is, to the existence of a 
faithful,
nondegenerate homomorphism of
$C_0(G^0)$ into the central multipliers $ZM(B)$. So, assume we have an
injective, nondegenerate homomorphism $\phi\colon C_0(G^0)\to ZM(B)$. For
$u\in G^0$ put
\begin{align*}
I_u&=\{f\in C_0(G^0):f(u)=0\}\\
K_u&=\phi(I_u)B&&\text{(a closed ideal of $B$)}\\
A_u&=B/K_u.
\end{align*}
Then put $\A=\bigcup_{u\in G^0}A_u$, and define $\Phi\colon
B\to\prod_{u\in G^0}A_u$ by
\[
\Phi(b)(u)=b+K_u.
\]
Then there is a unique topology on $\A$ making $\A$ an upper
semicontinuous $C^*$-bundle and each $\Phi(b)$ a continuous section, and
moreover $\Phi$ is an isomorphism of $B$ onto $\Gamma_0(\A)$. 

We need to relate the homomorphism $\phi$ to the action 
$(B,S,\beta)$. For $e\in E_S$ define the ideal
\[
C_e=\{f\in C_0(G^0):f=0\text{ off }e\}
\]
of $C_0(G^0)$. For our purposes, the appropriate connection 
between $\phi$ and $\beta$ is
\begin{equation}
\label{orthogonal}
\phi(C_e)B=B_e\righttext{for}e\in E_S,
\end{equation}
so we assume this henceforth.
Although we do not need it, we point out that \eqref{orthogonal} implies
$\phi$ is equivariant for $\beta$ and an obvious action of $S$ on
$C_0(G^0)$.

To simplify the writing, we use the isomorphism $\Phi\colon
B\to\Gamma_0(\A)$ to \emph{replace $B$ by $\Gamma_0(\A)$}. Then $\beta$
is an action of $S$ on $\Gamma_0(\A)$, and the homomorphism $\phi\colon
C_0(G^0)\to ZM(B)$ becomes the canonical embedding of $C_0(G^0)$ in
$ZM(\Gamma_0(\A))$. Since
\[
C_e\Gamma_0(\A)=\{f\in\Gamma_0(\A):f=0\text{ off }e\}
\righttext{for}e\in E_S,
\]
our hypothesis \eqref{orthogonal} tells us the ideals associated with the
inverse semigroup action $\beta$ are $B_e=\{f\in\Gamma_0(\A):f=0\text{
off }e\}$.

We want to construct an action $\alpha$ of the groupoid $G$ on the 
$C^*$-bundle $\A$. For a start, if $x\in G$ we need an isomorphism 
$\alpha_x$ of $A_{d(x)}$ onto $A_{r(x)}$. Take any $s\in S$ such that 
$x\in s$ (and 
note that such $s$ form a neighborhood base at $x$ 
in $G$).
For $u\in G^0$ we have
\[
K_u=\{f\in\Gamma_0(\A):f(u)=0\}.
\]
Furthermore, if $u\in e\in E_S$ we have
\[
B=B_e+K_u.
\]
Therefore, 
\[
A_u=B/K_u=(B_e+K_u)/K_u\cong B_e/(B_eK_u).
\]

\begin{lem}
With the above notation, there is a unique homomorphism
$\alpha_x$ 
from $A_{d(x)}$ to $A_{r(x)}$ such that
\[
\alpha_x(f(d(x)))=\beta_s(f)(r(x))
\righttext{for}x\in s\in S,f\in B_{d(s)}.
\]
\end{lem}

\begin{proof}
Identifying $A_{d(x)}$ with 
$B_{d(s)}/(B_{d(s)}K_{d(x)})$, 
and similarly for $A_{r(x)}$,
the conclusion of the lemma 
is equivalent to the assertion that there is a homomorphism 
$\alpha_x$ making the diagram
\[
\begin{diagram}
\node{ B_{d(s)}}
\arrow{e,t}{\beta_s}
\arrow{s}
\node{ B_{r(s)}}
\arrow{s}\\
\node{ B_{d(s)}\bigm/\bigl( B_{d(s)}K_{d(x)}\bigr)}
\arrow{e,b}{\alpha_x}
\node{ B_{r(s)}\bigm/\bigl( B_{r(s)}K_{r(x)}\bigr)}
\end{diagram}
\]
commute. For this we must show
\[
\beta_s\bigl( B_{d(s)}K_{d(x)}\bigr)
\subset B_{r(s)}K_{r(x)}.
\]
Take 
$f\in B_{d(s)}K_{d(x)}$. 
By density and continuity 
we can assume $\supp f\subset e$ for some $e\in E_S$ with $e\subset 
d(s)$ and $d(x)\notin e$. Then
\[
\beta_s(f)(r(x))=\beta_s\beta_e(f)(r(x))
=\beta_{se}(f)(r(x))=0,
\]
since
\[
\beta_{se}(f)\in B_{r(se)}= B_{ses^*}
\midtext{and}
r(x)=sd(x)s^*\notin ses^*.
\]
\end{proof}

\begin{thm}
Let $G$ be an $r$-discrete groupoid, let $S$ be 
a full inverse semigroup of open $G$-sets,
and let $\beta$ be an action 
of $S$ on a $C^*$-algebra $B$. Assume that there is an injective,
nondegenerate homomorphism $\phi$ of $C_0(G^0)$ into $ZM(B)$ such 
that $\phi(C_e)B=B_e$ for every $e\in E_S$. Then $B$ is 
isomorphic to $\Gamma_0(\A)$ for an upper semicontinuous 
$C^*$-bundle $\A$, and 
the map $\alpha$ defined in the above lemma is an action
of $G$ on $\A$.
\end{thm}

\begin{proof}
By the above discussion, the only thing left to check is that 
$\alpha$ is an action. We must check functoriality and continuity. 
Let $(x,y)\in G^2$ and $a\in A_{d(y)}$. Choose $s,t\in S$ such that 
$x\in s$ and $y\in t$, and then choose 
$f\in B_{d(st)}$ such 
that $f(d(y))=a$. Then $xy\in st$, so by the above 
lemma we have
\begin{align*}
\alpha_x\alpha_y(f(d(y)))
&=\alpha_x\bigl(\beta_t(f)(r(y))\bigr)
=\alpha_x\bigl(\beta_t(f)(d(x))\bigr)\\
&=\beta_s\beta_t(f)(r(x))
=\beta_{st}(f)(r(xy))\\
&=\alpha_{xy}\bigl(f(d(xy))\bigr)
=\alpha_{xy}\bigl(f(d(y))\bigr).
\end{align*}
Thus, $\alpha$ preserves compositions. To see that it preserves
identities, that is, $\alpha_u=\id_{A_u}$ for $u\in G^0$, just note that
$\alpha_u$ will be an idempotent surjection from $A_u$ to itself.

For the continuity, we appeal to \lemref{continuous}. Take $x\in G$, 
$a\in A_{d(x)}$, and $f,g\in\Gamma_c(\A)$ with $f(r(x))=\alpha_x(a)$ 
and $g(d(x))=a$. Cutting down $f$ and $g$, we can assume that 
$f\in B_{r(s)}$ and $g\in B_{d(s)}$ for some $s\in S$ 
with $x\in s$. Then for $y\in s$ we have
\begin{align*}
\norm{f(r(y))-\alpha_y(g(d(y)))}
&=\norm{f(r(y))-\beta_s(g)(r(y))}\\
&=\norm{\bigl(f-\beta_s(g)\bigr)(r(y))},
\end{align*}
which goes to $0$ as $y\to x$ since 
the norm is upper semicontinuous
and $f-\beta_s(g)$ is a continuous section which is $0$ at $r(x)$.
\end{proof}

\section{The isomorphism}
Let $G$ be an $r$-discrete groupoid, and let $S$ be 
a full inverse semigroup of open 
$G$-sets.
In the preceding two sections we established a correspondence between
the actions of $G$ and certain actions of $S$.  Our main result will be
that the associated crossed products are isomorphic.  To be specific, let
$\alpha$ be an action of $G$ on an upper semicontinuous $C^*$-bundle $\A$,
and let $\beta$ be the corresponding action of $S$ on the $C_0$-section
algebra $B:=\Gamma_0(\A)$. Recall that the ideals of $\beta$ are given by
\[
B_e:=\{f\in\Gamma_0(\A):f=0\text{ off }e\}\righttext{for}e\in E_S,
\]
and the partial automorphisms are given by
\[
\beta_s(f)(r(x))=\alpha_x(f(d(x)))
\righttext{for}x\in s\in S,f\in B_{d(s)}.
\]
We will show in \thmref{main} that $B\times_\beta S\cong \A\times_\alpha
G$.  
At the same time, we will fulfill our promise from \secref{gpoid}
by giving an independent proof that $\A\times_\alpha G$ exists, that is,
the $^*$-algebra $\Gamma_c(r^*\A)$ has an enveloping $C^*$-algebra. We
emphasize that our proof of this is completely independent of Renault's
(or anyone else's) decomposition theorem for representations of $G$. There
is no measure theory (other than the hypothesis that counting measure
is a Haar system on $G$); rather, the techniques are topological.  As a
(minor) byproduct, we have no separability requirements.

The crossed product $B\times_\beta S$ is the enveloping $C^*$-algebra of
a quotient of the finitely supported section algebra of the pull-back
bundle $r^*\B$. However, for our proof we will need to work with a
subbundle having incomplete fibers.
For $e\in E_S$ put $\Gamma_e(\A)=\{f\in\Gamma_c(\A):\supp f\subset e\}$,
and for
$s\in S$ put
\[
C_s=(\Gamma_{r(s)}(\A),s),
\]
giving a subbundle $\CC=\bigcup_{s\in S}C_s$ of $r^*\B$.

\begin{thm}
\label{bundle}
With the above notation, the
map $\Psi\colon\Gamma_c(\CC)\to\Gamma_c(r^*\A)$ defined on the
generators by
\[
\Psi(b,s)(x)=\begin{cases}b(r(x))&\text{if }x\in s,\\
0&\text{else}\end{cases}
\]
\textup(and extended additively\textup) 
is a surjective
$^*$-homomorphism with kernel $I_\beta\cap \Gamma_c(\CC)$.
\end{thm}

\begin{proof}
Since the range map $r$ takes each $G$-set $s\in S$ homeomorphically
onto $r(s)$, $\Psi$ takes each fiber $C_s$ of $\CC$
isometrically and isomorphically onto the linear subspace
\[
\Gamma_s(r^*\A):=\{f\in\Gamma_c(r^*\A):\supp f\subset s\}
\]
of $\Gamma_c(r^*\A)$. Since the elements of $S$ cover the groupoid $G$,
a standard partition of unity argument shows $\Psi$ maps
$\Gamma_c(\CC)$ onto $\Gamma_c(r^*\A)$.

Fix $(b,s),(c,t)\in \CC$. We have
\begin{align*}
\Psi\bigl((b,s)(c,t)\bigr)(x)
&=\Psi(\beta_s(\beta_{s^*}(b)c),st)(x)
\\&=\begin{cases}\beta_s(\beta_{s^*}(b)c)(r(x))&\text{if }x\in st,\\
0&\text{else.}\end{cases}
\end{align*}
Now, if $x\in st$ then $x$ factors \emph{uniquely} as $x=yz$ with $y\in
s$ and $z\in t$, and then
\begin{align*}
\beta_s(\beta_{s^*}(b)c)(r(x))
&=\beta_s(\beta_{s^*}(b)c)(r(y))
=\alpha_y\bigl((\beta_{s^*}(b)c)(d(y))\bigr)
\\&=\alpha_y\bigl(\beta_{s^*}(b)(d(y))c(d(y))\bigr)
\\&=\alpha_y\bigl(\alpha_{y^{-1}}(b(r(y)))\bigr)\alpha_y\bigl(c(d(y))\bigr)
\\&=b(r(y))\alpha_y\bigl(c(r(y^{-1}x))\bigr)
\\&=\Psi(b,s)(y)\alpha_y\bigl(\Psi(c,t)(y^{-1}x)\bigr)
\\&=\sum_{r(w)=r(x)}\Psi(b,s)(w)\alpha_w\bigl(\Psi(c,t)(w^{-1}x)\bigr)
\\&=\bigl(\Psi(b,s)\Psi(c,t)\bigr)(x).
\end{align*}
On the other hand,
\[
\supp\bigl(\Psi(b,s)\Psi(c,t)\bigr)
\subset\bigl(\supp\Psi(b,s)\bigr)\bigl(\supp\Psi(c,t)\bigr)
\subset st,
\]
so if $x\notin st$ then
\[
0=\bigl(\Psi(b,s)\Psi(c,t)\bigr)(x).
\]
Hence, $\Psi$ is multiplicative. A similar computation shows $\Psi$
preserves adjoints, so $\Psi$ is a \mbox{$^*$-ho}\-mo\-mor\-phism.

It remains to show the kernel of the map
$\Psi\colon\Gamma_c(\CC)\to\Gamma_c(r^*\A)$ is the ideal
$I_\beta\cap\Gamma_c(\CC)$.  This
will
involve a couple of mildly fussy partition-of-unity arguments, so we have
made an attempt to isolate the hard bit by factoring the map $\Psi$
through an auxiliary bundle: 
let $\X$ be the 
bundle over $S$ with (incomplete) fibers
$\{(\Gamma_s(r^*\A),s)\}_{s\in S}$. Then 
define
$\Theta\colon\Gamma_c(\CC)\to\Gamma_c(\X)$ and
$\Lambda\colon\Gamma_c(\X)\to\Gamma_c(r^*\A)$ by
\[
\Theta(b,s)=(\Psi(b),s)\midtext{and}\Lambda(f,s)=f,
\]
so that $\Psi=\Lambda\circ\Theta$. The map $\Theta$ is a linear
isomorphism of $\Gamma_c(\CC)$ onto $\Gamma_c(\X)$, since it is induced
by the identity map on the base space $S$ of the bundles $\CC$ and $\X$,
and by linear isomorphisms between the fibers $(\Gamma_{r(s)}(\A),s)$ and
$(\Gamma_s(r^*\A),s)$. Moreover, $\Theta$ takes $I_\beta\cap\Gamma_c(\CC)$
onto the span of $\{(f,s)-(f,t):f\in\Gamma_s(r^*\A),s\le t\}$. Thus,
it remains to show the kernel of $\Lambda$ coincides with this span.

We first show
\begin{equation}
\label{step 1}
\ker\Lambda=\spn\{(f,s)-(f,t):f\in\Gamma_{s\cap t}(r^*\A)\}.
\end{equation}
Let $I$ denote the right hand side. Certainly
$I\subset\ker\Lambda$.  For the opposite inclusion, suppose
$\Lambda\bigl(\sum_1^n(f_i,s_i)\bigr)=0$. We need to show
$\sum_1^n(f_i,s_i)\in I$. We have $\sum_1^nf_i=0$, so if $n=1$ then
$f_1=0$ and so $\sum_1^n(f_i,s_i)=(f_1,s_1)=0$. Hence, we can assume
$n>1$. Let $\Omega$ denote the family of subsets of $\{1,\dots,n\}$
with cardinality at least $2$, and for $\omega\in\Omega$ put
\[
V_\omega=\biggl(\bigcap_{i\in\omega}s_i\biggr)\biggm\backslash
\biggl(\bigcup_{i\notin\omega}\supp f_i\biggr).
\]
Then $\{V_\omega\}_{\omega\in\Omega}$ forms an open cover of
$\bigcup_1^n\supp f_i$, since if $f_i(x)\not=0$ then
also $f_j(x)\not=0$ for at least one $j\not=i$, and upon taking limits
we get $\supp f_i\subset\bigcup_{j\not=i}\supp f_j$, so if
$x\in\supp f_i$ and $\omega=\{j:x\in s_j\}$, then $\omega\in\Omega$ and
$x\in V_\omega$.
Choose a partition
of unity $\{\phi_\omega\}_{\omega\in\Omega}$ subordinate to the open
cover $\{V_\omega\}_{\omega\in\Omega}$ of $\bigcup_1^n\supp f_i$.
Note that whenever $i\notin\omega$ we have
$f_i\phi_\omega=0$ since $\supp f_i\cap V_\omega=\emptyset$. Hence,
\[
\sum_{i\in\omega}f_i\phi_\omega
=\sum_{i=1}^nf_i\phi_\omega
=(0)\phi_\omega
=0\righttext{for all}\omega\in\Omega.
\]
We have
\begin{align*}
\sum_1^n(f_i,s_i)
&=\sum_1^n\biggl(\sum_{\omega\in\Omega}f_i\phi_\omega,s_i\biggr)
=\sum_\omega\sum_1^n(f_i\phi_\omega,s_i)
\\&=\sum_\omega\sum_{i\in\omega}(f_i\phi_\omega,s_i).
\end{align*}
Fix $\omega\in\Omega$, and pick any two distinct elements $j,k$ of
$\omega$. Then
\begin{align*}
\sum_{i\in\omega}(f_i\phi_\omega,s_i)
&=(f_j\phi_\omega,s_j)+(f_k\phi_\omega,s_k)
+\sum_{i\in\omega\setminus\{j,k\}}(f_i\phi_\omega,s_i)
\\&=(f_j\phi_\omega,s_j)+(f_k\phi_\omega,s_k)
+\sum_{i\in\omega\setminus\{j,k\}}(f_i\phi_\omega,s_i)
\\&\quad+\sum_{i\in\omega\setminus\{j,k\}}(f_i\phi_\omega,s_k)
-\sum_{i\in\omega\setminus\{j,k\}}(f_i\phi_\omega,s_k)
\\&=(f_j\phi_\omega,s_j)
+\sum_{i\in\omega\setminus\{j\}}(f_i\phi_\omega,s_k)
\\&\quad+\sum_{i\in\omega\setminus\{j,k\}}
\bigl((f_i\phi_\omega,s_i)-(f_i\phi_\omega,s_k)\bigr)
\\&=(f_j\phi_\omega,s_j)-(f_j\phi_\omega,s_k)
\\&\quad+\sum_{i\in\omega\setminus\{j,k\}}
\bigl((f_i\phi_\omega,s_i)-(f_i\phi_\omega,s_k)\bigr),
\end{align*}
because $f_j\phi_\omega+\sum_{i\in\omega\setminus\{j\}}f_i\phi_\omega=0$.
Since $\supp f_i\phi_\omega\subset s_l$ for every $i,l\in\omega$, we
conclude that $\sum_{i\in\omega}(f_i\phi_\omega,s_i)$ is an element of
$I$, so we have shown \eqref{step 1}.

Now put
\[
J=\spn\{(f,s)-(f,t):f\in\Gamma_s(r^*\A),s\le t\}.
\]
Clearly $J\subset\ker\Lambda$. For the opposite containment, by the
above argument it suffices to show that if $f\in\Gamma_{s\cap
t}(r^*\A)$ then $(f,s)-(f,t)\in J$. Since $S$ is a base for the
topology of $G$, we can find $s_1,\dots,s_n\in S$ such that
\[
\supp f\subset\bigcup_1^ns_i\subset s\cap t.
\]
choose a partition of unity $\{\phi_i\}_1^n$ subordinate to the open
cover $\{s_i\}_1^n$ of $\supp f$. We have
\begin{align*}
(f,s)-(f,t)
&=\sum_1^n(f\phi_i,s)-\sum_1^n(f\phi_i,t)
\\&=\sum_1^n(f\phi_i,s)-\sum_1^n(f\phi_i,s_i)
+\sum_1^n(f\phi_i,s_i)-\sum_1^n(f\phi_i,t)
\\&=\sum_1^n\bigl((f\phi_i,s)-(f\phi_i,s_i)\bigr)
+\sum_1^n\bigl((f\phi_i,s_i)-(f\phi_i,t)\bigr).
\end{align*}
Since each $\phi_i$ has support contained in $s_i\cap s\cap t$, the
latter sums are elements of $J$, and we are done.
\end{proof}

\begin{thm}
\label{main}
Let $\alpha$ be an action of an $r$-discrete groupoid $G$ on
an upper semicontinuous $C^*$-bundle $\A$, $S$ a full inverse
semigroup of open $G$-sets \textup(as in \defnref{full}\textup), and
$\beta$ the associated action of $S$ on $B:=\Gamma_0(\A)$ \textup(as in
\thmref{G-S}\textup).
Then the $^*$-algebra
$\Gamma_c(r^*\A)$ has an enveloping $C^*$-algebra $\A\times_\alpha
G$. Moreover, the map $\Psi$ of \thmref{bundle} extends uniquely to an
isomorphism of $B\times_\beta S$ onto $\A\times_\alpha G$.
\end{thm}

\begin{proof}
By \thmref{bundle}, the map $\Psi$ factors through
an isomorphism $\Psi'$ of the quotient $\Gamma_c(\CC)/(I_\beta\cap
\Gamma_c(\CC))$ onto $\Gamma_c(r^*\A)$.  Since each $\Gamma_e(\A)$ is a
dense ideal of $B_e$ and $\beta_s(\Gamma_{d(s)}(\A))=\Gamma_{r(s)}(\A)$
for every $s\in S$, \corref{S envelope} tells us $B\times_\beta S$
is the enveloping $C^*$-algebra of $\Gamma_c(\CC)/(I_\beta\cap
\Gamma_c(\CC))$.  The result follows.
\end{proof}

\section{Application} We show how 
\thmref{main} allows us to recover 
Paterson's representation \cite{pat:book} 
of the $C^{*}$-algebra of an $r$-discrete groupoid as a 
semigroup crossed product.  Since our groupoids are 
Hausdorff and Paterson requires only the unit space of 
the groupoids to be Hausdorff, we cannot get his 
theorem in full generality.  We believe the connection 
between groupoid and inverse semigroup crossed products 
should also work for non-Hausdorff groupoids.  

If $G$ is 
a not necessarily Hausdorff
$r$-discrete groupoid then Paterson \cite{pat:book}
calls an inverse subsemigroup $S$ of $G^{op}$ 
{\em additive\/} if $S$ is a base 
for the topology of $G$ and $s$, 
$t\in S$ with $s\cup t\in G^{op}$ implies 
$s\cup t\in S$.  
Note that additivity is a strictly stronger condition than fullness in
the sense of \defnref{full}.
Paterson shows 
\cite[Theorem 3.3.1]{pat:book}
that if $S$ is an additive inverse subsemigroup of $G^{op}$ 
then $C^{*}(G)$ 
is isomorphic to 
$C_0(G^0)\times_\beta S$ 
if $S$ acts on $C_0(G^0)$ canonically (see below).  
We can deduce the same result if we assume that $G$ is also 
Hausdorff, and in fact we can get away with slightly less than 
additivity:

\begin{thm}
Let $G$ be an $r$-discrete 
\textup(Hausdorff\textup)
groupoid and let $S$ be a full inverse semigroup of open 
$G$-sets.
Then $S$ has an action $\beta$ on $B:=C_0(G_0)$ defined by  
\[\beta_s(f)(u)=\begin{cases}f(s^*us)&\text{if }u\in 
r(s),\\
     0&\text{else},\end{cases}\]
for 
$f\in B_{d(s)}:=\{f\in C_0(G^0):f=0\text{ off } d(s)\}$,
and $C^{*}(G)$ is isomorphic to $B\times_{\beta}S$.
\end{thm}

\begin{proof}
$G$ has an action $\alpha$ on the trivial $C^{*}$-bundle 
$\A={\C}\times G^0$, 
where $\alpha_x:A_{d(x)}\to A_{r(x)}$ is the identity map between two 
copies of ${\mathbb C}$.  It is clear that $C^{*}(G)$ is isomorphic to 
$\A\times_{\alpha}G$.  
By \thmref{main}
$\A\times_{\alpha}G$ is isomorphic to $B\times_{
\beta}S$.
Since $\Gamma_0(\A)$ is isomorphic to $C_0(G^0)$ and $\alpha_
x$ is 
the identity map for all $x\in G$, $\beta$ is exactly the canonical 
action of $S$ on $C_0(G^0)$ used by Paterson.  
\end{proof}


\providecommand{\bysame}{\leavevmode\hbox to3em{\hrulefill}\thinspace}

\end{document}